\documentclass[a4paper,reqno]{amsart}

\usepackage{amsmath}
\usepackage{amsfonts}
\usepackage{a4wide}
\usepackage{amssymb}
\usepackage{amsthm}

\theoremstyle{plain}

\newtheorem{teo}{Theorem}[section]

\newtheorem{cor}[teo]{Corollary}
\newtheorem{ackn}{Acknowledgments\!}

\theoremstyle{definition}

\theoremstyle{remark}

\newtheorem{rem}[teo]{Remark}

\numberwithin{equation}{section}

\def\SS{{{\mathbb S}}}

\def\RR{{\mathbb R}}

\def\gt{\widetilde{g}}

\def\RRR{{\mathrm R}}
\def\WWW{{\mathrm W}}

\def\HHH{{\mathrm H}}

\def\Ric{{\mathrm {Ric}}}

\def\SSS{{\mathrm S}}
\def\CCC{{\mathrm C}}
\def\Rm{{\mathrm {Rm}}}

\newcommand{\ti}[1]{\widetilde{#1}}

\title[Generalized Quasi--Einstein Manifolds with Harmonic Weyl Tensor]{Generalized Quasi--Einstein Manifolds \\with Harmonic Weyl Tensor}
\date{\today}

\author[Giovanni Catino]{Giovanni Catino}
\address[Giovanni Catino]{SISSA -- International School for Advanced Studies, Via Bonomea 265 , Trieste, Italy, 34136}
\email[]{catino@sissa.it}

\date{\today}

\begin{document}

\begin{abstract} In this paper we introduce the notion of generalized quasi--Einstein manifold, that generalizes the concepts of Ricci soliton, Ricci almost soliton and quasi--Einstein manifolds. We prove that a complete generalized quasi--Einstein manifold with harmonic Weyl tensor and with zero radial Weyl curvature, is locally a warped product with $(n-1)$--dimensional Einstein fibers. In particular, this implies a local characterization for locally conformally flat gradient Ricci almost solitons, similar to that proved for gradient Ricci solitons.
\end{abstract}

\maketitle

\medskip

\section{Introduction}

\medskip

In recent years, much attention has been given to the classification of Riemannian manifolds admitting an Einstein--like structure. In this paper we will define a class of Riemannian metrics which naturally generalizes the Einstein condition. More precisely, we say that a complete Riemannian manifold $(M^{n},g)$, $n\geq 3$,  is a {\em generalized quasi--Einstein manifold}, if there exist three smooth functions $f,\mu,\lambda$ on $M$, such that 
\begin{equation}\label{gqe}
\Ric + \nabla^{2} f - \mu\, df \otimes df = \lambda g \,.
\end{equation}
Natural examples of GQE manifolds are given by Einstein manifolds (when $f$ and $\lambda$ are two constants), gradient Ricci solitons (when $\lambda$ is constant and $\mu=0$), gradient Ricci almost solitons (when $\mu=0$, see~\cite{pigrimsett1}) and quasi--Einstein manifolds (when $\mu$ and $\lambda$ are two constants, see~\cite{caseshuwei}~\cite{mancatmazz1}~\cite{HePetWylie}). We will call a GQE manifolds {\em trivial}, if the function $f$ is constant. This will clearly imply that $g$ is an Einstein metric.

\medskip

The Riemann curvature operator of a Riemannian manifold $(M^n,g)$ is defined 
as in~\cite{gahula} by
$$
\mathrm{Riem}(X,Y)Z=\nabla_{Y}\nabla_{X}Z-\nabla_{X}\nabla_{Y}Z+\nabla_{[X,Y]}Z\,.
$$ 
In a local coordinate system the components of the $(3,1)$--Riemann 
curvature tensor are given by
$\RRR^{d}_{abc}\tfrac{\partial}{\partial
  x^{d}}=\mathrm{Riem}\big(\tfrac{\partial}{\partial
  x^{a}},\tfrac{\partial}{\partial
  x^{b}}\big)\tfrac{\partial}{\partial x^{c}}$ and we denote by
$\RRR_{abcd}=g_{de}\RRR^{e}_{abc}$ its $(4,0)$--version.

\medskip

{\em In all the paper the Einstein convention of summing over the repeated indices will be adopted.}

\medskip

With this choice, for the sphere $\SS^n$ we have
${\mathrm{Riem}}(v,w,v,w)=\RRR_{abcd}v^aw^bv^cw^d>0$. The Ricci tensor is obtained by the contraction 
$\RRR_{ac}=g^{bd}\RRR_{abcd}$ and $\RRR=g^{ac}\RRR_{ac}$ will 
denote the scalar curvature. The so called Weyl tensor is then 
defined by the following decomposition formula (see~\cite[Chapter~3,
Section~K]{gahula}) in dimension $n\geq 3$,
\begin{eqnarray*}
\WWW_{abcd}=&\,\RRR_{abcd}+\frac{\RRR}{(n-1)(n-2)}(g_{ac}g_{bd}-g_{ad}g_{bc})
- \frac{1}{n-2}(\RRR_{ac}g_{bd}-\RRR_{ad}g_{bc}
+\RRR_{bd}g_{ac}-\RRR_{bc}g_{ad})\,.
\end{eqnarray*}

We recall that a Riemannian metric has {\em harmonic Weyl tensor} if the divergence of $\WWW$ vanishes. In dimension three this condition is equivalent to local conformally flatness. Nevertheless, when $n\geq 4$, harmonic Weyl tensor is a weaker condition since locally conformally flatness is equivalent to the vanishing of the Weyl tensor. 

\medskip

In this paper we will give a local characterization of generalized quasi--Einstein manifolds with harmonic Weyl tensor and such that $\WWW(\nabla f,\cdot, \cdot, \cdot)=0$. As we have seen, this class includes the case of locally conformally flat manifolds.

\begin{teo}\label{main}
Let $(M^{n},g)$, $n\geq 3$, be a generalized quasi--Einstein manifold with harmonic Weyl tensor and $\WWW(\nabla f,\cdot, \cdot, \cdot)=0$. Then, around any regular point of $f$, the manifold $(M^{n},g)$ is locally a warped product with $(n-1)$--dimensional Einstein fibers.
\end{teo}

\begin{rem} We notice that the hypothesis $\WWW(\nabla f,\cdot, \cdot, \cdot)=0$ cannot be removed. Indeed, if we consider the gradient shrinking soliton on $M=\RR^{k}\times \mathbb{S}^{n-k}$, for $n\geq 4$ and $k\geq 2$, defined by the product metric $g=dx^{1}\otimes\cdots\otimes dx^{k}+g_{\mathbb{S}^{n-k}}$ and the potential function
$$
f \,=\, \tfrac{1}{2}\,\big( \,|x^{1}|^{2}+\dots|x^{k}|^{2} \,\big) \,,
$$
it is easy to verify that $(M^{n},g)$ has harmonic Weyl tensor, since it is the product of two Einstein metrics, whereas the radial part of the Weyl tensor $\WWW(\nabla f,\cdot, \cdot, \cdot)$ does not vanish.
\end{rem}

\begin{rem} Theorem~\ref{main} generalizes the results obtained for gradient Ricci solitons (see~\cite{caochen} and~\cite{mancat1}) and, recently, for quasi--Einstein manifolds (see~\cite{mancatmazz1}).
\end{rem}

As an immediate corollary, we have that a locally conformally flat generalized quasi--Einstein manifold is, locally, a warped product with $(n-1)$--dimensional fibers of constant sectional curvature. In particular, we can prove a local characterization for locally conformally flat Ricci almost solitons (which have been introduced in~\cite{pigrimsett1}), similar to the one for Ricci solitons (\cite{caochen}~\cite{mancat1}). 

\begin{cor} 
Let $(M^{n},g)$, $n\geq 3$, be a locally conformally flat gradient Ricci almost soliton. Then, around any regular point of $f$, the manifold $(M^{n},g)$ is locally a warped product with $(n-1)$--dimensional fibers of constant sectional curvature.
\end{cor}

If $n=4$, since a three dimensional Einstein manifold has constant sectional curvature, we get the following 

\begin{cor}
Let $(M^{4},g)$, be a four dimensional generalized quasi--Einstein manifold with harmonic Weyl tensor and $\WWW(\nabla f,\cdot, \cdot, \cdot)=0$. Then, around any regular point of $f$, the manifold $(M^{4},g)$ is locally a warped product with three dimensional fibers of constant sectional curvature. In particular, if it is nontrivial, then $(M^{4},g)$ is locally conformally flat.
\end{cor}

Now, using the classification of locally conformally flat gradient steady Ricci solitons (see again~\cite{caochen} and \cite{mancat1}), we obtain

\begin{cor}
Let $(M^{4},g)$, be a four dimensional gradient steady Ricci soliton with harmonic Weyl tensor and $\WWW(\nabla f,\cdot, \cdot, \cdot)=0$. Then $(M^{4},g)$ is either Ricci flat or isometric to the Bryant soliton. 
\end{cor}

\medskip

\section{Proof of Theorem~\ref{main}}

\medskip

Let $(M^{n},g)$, $n\geq 3$, be a generalized quasi--Einstein manifold with harmonic Weyl tensor and satisfying $\WWW(\nabla f,\cdot, \cdot, \cdot)=0$. 
If $n=3$, we have that $g$ is locally conformally flat, while if $n\geq 4$, one has
\begin{align*}
0=&\,\nabla^d\WWW_{abcd}\\
=&\,\nabla^d\Bigl(\RRR_{abcd}+\frac{\RRR}{(n-1)(n-2)}(g_{ac}g_{bd}-g_{ad}g_{bc})
- \frac{1}{n-2}(\RRR_{ac}g_{bd}-\RRR_{ad}g_{bc}
+\RRR_{bd}g_{ac}-\RRR_{bc}g_{ad})\Bigr)\\
=&\,-\nabla_a\RRR_{bc}+\nabla_b\RRR_{ac}
+\frac{\nabla_b\RRR}{(n-1)(n-2)}g_{ac}
-\frac{\nabla_a\RRR}{(n-1)(n-2)}g_{bc}\\
&\,- \frac{1}{n-2}(\nabla_b\RRR_{ac}-\nabla^d\RRR_{ad}g_{bc}
+\nabla^d\RRR_{bd}g_{ac}-\nabla_a\RRR_{bc}g_{ad})\\
=&\,-\frac{n-3}{n-2}(\nabla_a\RRR_{bc}-\nabla_b\RRR_{ac})
+\frac{\nabla_b\RRR}{(n-1)(n-2)}g_{ac}
-\frac{\nabla_a\RRR}{(n-1)(n-2)}g_{bc}\\
&\,+\frac{1}{2(n-2)}(\nabla_a\RRR g_{bc}/2
-\nabla_b\RRR g_{ac}/2)\\
=&\,-\frac{n-3}{n-2}\Bigl[
\nabla_a\RRR_{bc}-\nabla_b\RRR_{ac}-\frac{(\nabla_a\RRR g_{bc}-\nabla_b\RRR g_{ac})}{2(n-1)}\Bigr]\\
=&\,-\frac{n-3}{n-2}\CCC_{cba}\\
=&\,-\frac{n-3}{n-2}\CCC_{abc} \,,
\end{align*}
where $\CCC$ is the Cotton tensor
$$
\CCC_{abc} \,=\, \nabla_{c} \RRR_{ab} - \nabla_{b} \RRR_{ac} -
\tfrac{1}{2(n-1)} \big( \nabla_{c} \RRR \, g_{ab} - \nabla_{b} \RRR \,
g_{ac} \big)\,.
$$
Hence, if $n\geq 3$, harmonic Weyl tensor is equivalent to the vanishing of the Cotton tensor. 

Now, the condition $\WWW(\nabla f,\cdot, \cdot, \cdot)=0$ implies that the conformal metric 
$$
\gt \,= \,e^{-\frac{2}{n-2}f}g
$$ 
has harmonic Weyl tensor. Indeed, from the conformal transformation law for the Cotton tensor (see Appendix), one has that, if $n\geq 4$, then
$$
(n-2)\,\widetilde{\CCC}_{abc} \, = \, (n-2)\,\CCC_{abc}+\tfrac{1}{n-2}\WWW_{abcd}\nabla^{d}f = 0 \,,
$$ 
whereas $\widetilde{\CCC}_{abc}=\CCC_{abc}=0$ in three dimensions. Hence, from the definition of the Cotton tensor, we can observe that the Schouten tensor of $\gt$ defined by
$$
\SSS_{\gt} \, = \, \tfrac{1}{n-2}\big(\,\Ric_{\gt}-\tfrac{1}{2(n-1)}\,\RRR_{\gt}\,\,\gt \,\big)
$$ 
is a Codazzi tensor, i.e. it satisfies the equation
$$
(\nabla_{X}\SSS)\,Y \,=\, (\nabla_{Y}\SSS)\,X\,, \quad\,\hbox{for all}\,\, X,Y\in TM\,.
$$
(see~\cite[Chapter~16, Section~C]{besse} for a general overview on Codazzi
tensors). 

Moreover, from the structural equation of generalized quasi--Einstein manifolds~\eqref{gqe}, the expression of the Ricci tensor of the conformal metric $\gt$ takes the form
\begin{eqnarray*}
\Ric_{\widetilde{g}}&=&\Ric_{g} +\nabla^{2} f +\tfrac{1}{n-2} df \otimes
df +\tfrac{1}{n-2} \big(\Delta f - |\nabla f|^{2} \big) g\\
&=& \, \big(\mu + \tfrac{1}{n-2}\big) df \otimes df + \tfrac{1}{n-2}\big(\Delta f - |\nabla f|^{2}+(n-2)\lambda\big) \, e^{\frac{2}{n-2}f}\,\gt\,.
\end{eqnarray*}
Then, at every regular point $p$ of $f$, the Ricci tensor of $\gt$ either has a unique eigenvalue or has two distinct eigenvalues $\eta_{1}$ and $\eta_{2}$ of multiplicity $1$ and $(n-1)$ respectively. In both cases, $\nabla f/|\nabla f|_{\gt}$ is an eigenvector of the Ricci tensor of $\gt$. For every point in $\Omega=\{p\in M \,|\, p \,\,\hbox{regular point}, \eta_{1}(p)\neq \eta_{2}(p)\}$ also the Schouten tensor $\SSS_{\gt}$ has two distinct eigenvalues $\sigma_{1}$ of multiplicity one and $\sigma_{2}$ of multiplicity $(n-1)$, with same eigenspaces of $\eta_{1}$ and $\eta_{2}$ respectively.
Splitting results for Riemannian manifolds admitting a Codazzi tensor with only two distinct eigenvalues were obtained by Derdzinski~\cite{derdz3} and Hiepko--Reckziegel~\cite{hiepkoreck} (see again~\cite[Chapter~16, Section~C]{besse} for further discussion). 

From Proposition 16.11 in~\cite{besse} (see also~\cite{derdz3}) we know that the tangent bundle of a neighborhood of $p$ splits as the orthogonal direct sum of two integrable eigendistributions, a line field $V_{\sigma_{1}}$, and a codimension one distribution $V_{\sigma_{2}}$ with totally umbilic leaves, in the sense that the second fundamental form $\ti{h}$ of each leaves is proportional to the metric $\gt$ (with abuse of notation, we will call $\gt$ also the induced metric on the leaves of $V_{\sigma_{2}}$). We will denote by $\ti{\nabla}$ the Levi--Civita connection of the metric $\gt$ on $M$ and by $\ti{\nabla}^{\sigma_{2}}$ the induced Levi--Civita connection of the induced metric $\gt$ on the leaves of $V_{\sigma_{2}}$. In a suitable local chart $x^{1},x^{2},\dots,x^{n}$ with $\partial/\partial x^{1} \in V_{\sigma_{1}}$, $\partial/\partial x^{i} \in V_{\sigma_{2}}$ (in the sequel $i,j,k$ will range over $2,\dots,n$), we have $\gt_{1i}=0$. Since $V_{\sigma_{2}}$ is totally umbilic, we have
\begin{equation}\label{umbilic}
\ti{h}_{ij} \,= \,-\big\langle \ti{\nabla}^{\sigma_{2}}_{\tfrac{\partial}{\partial x^{i}}} \tfrac{\partial}{\partial x^{j}}, \tfrac{\partial}{\partial x^{1}} \big\rangle = -\ti{\Gamma}^{1}_{ij} \,\,\gt_{11} = \tfrac{\ti{\HHH}}{n-1} \ti{g}_{ij} \,,
\end{equation}
where $\ti{\HHH}$ will denote the mean curvature function. We recall that, from the Codazzi--Mainardi equation (see Theorem 1.72 in~\cite{besse}), one has
\begin{equation}\label{codman}
\big(\ti{\nabla}^{\sigma_{2}}_{\tfrac{\partial}{\partial x^{i}}}\ti{h}\big) \big(\tfrac{\partial}{\partial x^{j}},\tfrac{\partial}{\partial x^{k}}\big) - \big(\ti{\nabla}^{\sigma_{2}}_{\tfrac{\partial}{\partial x^{j}}}\ti{h}\big) \big(\tfrac{\partial}{\partial x^{i}},\tfrac{\partial}{\partial x^{k}}\big) \, = \,  \big\langle \,\ti{\Rm} \big(\tfrac{\partial}{\partial x^{i}},\tfrac{\partial}{\partial x^{j}}\big) \tfrac{\partial}{\partial x^{k}}, \tfrac{\partial}{\partial x^{1}} \, \big\rangle \,.  
\end{equation}
On the other hand, tracing with the metric $\gt$, and using the umbilic property~\eqref{umbilic}, we get
\begin{eqnarray*}
\big(\ti{\nabla}^{\sigma_{2}}_{\tfrac{\partial}{\partial x^{i}}}\ti{h}\big) \big(\tfrac{\partial}{\partial x^{j}},\tfrac{\partial}{\partial x^{i}}\big) - \big(\ti{\nabla}^{\sigma_{2}}_{\tfrac{\partial}{\partial x^{j}}}\ti{h}\big) \big(\tfrac{\partial}{\partial x^{i}},\tfrac{\partial}{\partial x^{i}})  \,=\,  \tfrac{1}{n-1}\,\partial_{j}\ti{\HHH} - \partial_{j}\ti{\HHH} \,=\, \tfrac{2-n}{n-1} \, \partial_{j}\ti{\HHH} \,.
\end{eqnarray*}
Using equation~\eqref{codman}, we obtain
$$
\tfrac{2-n}{n-1} \, \partial_{j}\ti{\HHH} \,=\, \Ric_{\gt} \big(\tfrac{\partial}{\partial x^{j}},\tfrac{\partial}{\partial x^{1}} \big) \,=\,  0 \,,
$$
which implies that the mean curvature $\ti{\HHH}$ is constant on each leaves of $V_{\sigma_{2}}$. Now, from Proposition 16.11 (ii) in~\cite{besse}, one has that
$$
\ti{\HHH} = \tfrac{1}{\sigma_{1}-\sigma_{2}}\, \partial_{1} \,\sigma_{2} \,.
$$
The facts that both $\ti{\HHH}$ and $\sigma_{2}$ are constant on each leaves of $V_{\sigma_{2}}$ imply that $\partial_{j} \, \sigma_{1}=0$, for every $j=2,\dots,n$. This is equivalent to say that $V_{\sigma_{1}}$ has to be a geodesic line distribution, which clearly implies $\Gamma^{j}_{00}=0$, i.e. $\partial_{j} \,g_{11} =0$. Equation~\eqref{umbilic} yields 
$$
\partial_{1} \gt_{ij} \,=\, -2 \,\ti{\Gamma}^{1}_{ij} = 2 \,\gt_{11}^{-1}\, \tfrac{\ti{\HHH}}{n-1}\, \gt_{ij} \,.
$$
Since $\ti{\HHH}$ and $g_{11}$ are constant along $V_{\sigma_{2}}$, one has
$$ 
\partial_{1} \gt_{ij}(x^{1},\dots,x^{n}) \,=\, \varphi(x^{1}) \,\gt_{ij}(x^{1},\dots,x^{n})
$$ for some function $\varphi$ depending only on the $x^{1}$ variable. Choosing a function $\psi=\psi(x^{1})$, such that $\tfrac{d\,\psi}{dx^{1}} = \varphi$, we have $\partial_{1} (e^{-\psi}\,\gt_{ij})=0$, which means that 
$$
\gt_{ij}(x^{1},\dots,x^{n})\,=\,e^{\psi(x^{1})} \, G_{ij}(x^{2},\dots,x^{n}) \,,
$$
for some $G_{ij}$. This implies that the manifold $(M^{n},\gt)$, locally around every regular point of $f$, has a warped product representation with $(n-1)$--dimensional fibers. By the structure of the conformal deformation, this conclusion also holds for the original Riemannian manifold $(M^{n}, g)$. Now, the fact that $g$ has harmonic Weyl tensor, implies that the $(n-1)$--dimensional fibers  are Einstein manifolds (there are a lot of papers where this computation is done, for instance see~\cite{geba}). 

This completes the proof of Theorem~\eqref{main}.

\medskip

\section*{Appendix}\label{appendix}

\medskip

\setcounter{section}{1}
\setcounter{teo}{0}

\noindent{\bf Lemma.} \,{\em The Cotton tensor $\CCC_{abc}$ is pointwise conformally invariant in dimension three, whereas if $n\geq 4$, for $\ti{g}=e^{-2u}g$, we have
$$
(n-2)\,\ti{\CCC}_{abc} \,= \, (n-2)\,\CCC_{abc}+\WWW_{abcd}\nabla^{d}u \,.
$$}
\vspace{-0.5cm}
\begin{proof}
The proof is a straightforward computation. Let $\ti{g}=e^{-2u}\,g$, then for the Schouten tensor $\SSS = \frac{1}{n-2} \big(\Ric - \frac{1}{2(n-1)}\RRR \,g \big)$ we
have the conformal transformation rule
\begin{equation}\label{schconf}
\ti{\SSS} \,=\, \SSS + \nabla^{2} u + du \otimes du - \tfrac{1}{2} |\nabla
u|^{2} g \,.
\end{equation}
The Cotton tensor of the metric $\ti{g}$ is defined by
$$
(n-2)\, \ti{\CCC}_{abc} \,=\, \ti{\nabla}_{c} \ti{\SSS}_{ab} - \ti{\nabla}_{b}
\ti{\SSS}_{ac} \,.
$$
Moreover one can see that
\begin{eqnarray*}
\ti{\nabla}_{c} \ti{\SSS}_{ab} 
&=& \nabla_{c} \SSS_{ab} + \nabla_{c} \nabla_{a} \nabla_{b} u +\nabla_{c}
\nabla_{a} u \, \nabla_{b} u + \nabla_{c} \nabla_{b} u \, \nabla_{a} u
- \nabla_{c} \nabla_{d} u\, \nabla_{d} u \, g_{ab}+\\
&& + \, \ti{\SSS}_{bc} \nabla_{a} u + \ti{\SSS}_{ac} \nabla_{b} u + 
\ti{\SSS}_{ab} \nabla_{c} u - \ti{\SSS}_{bd} \nabla_{d} u \, g_{ac} -
\ti{\SSS}_{ad} \nabla_{d} u \, g_{bc} \,.
\end{eqnarray*}
Computing in the same way the term $\ti{\nabla}_{b} \ti{\SSS}_{ac}$,
substituting in the previous formula $\ti{\SSS}$ with~\eqref{schconf} and
using the fact that
\begin{eqnarray*}
\nabla_{c} \nabla_{b} \nabla_{a} u - \nabla_{b} \nabla_{c} \nabla_{a}
u &=& \RRR_{cbad} \nabla^{d} u \,\,\, = \,\,\, \RRR_{abcd} \nabla^{d} u \\
&=& \WWW_{abcd} \nabla^{d} u + \SSS_{ac} \nabla_{b} u - \SSS_{cd} \nabla_{d} u
\, g_{ab} +\SSS_{bd} \nabla_{d} u \, g_{ac} - \SSS_{ab} \nabla_{c} u \,,
\end{eqnarray*}
(we recall that $\WWW$ is zero in dimension three) one obtains the result.
\end{proof}

\bigskip

\begin{ackn} The author is partially supported by the Italian project FIRB--IDEAS ``Analysis and Beyond''. The author wishes to thank Sara Rizzi for helpful remarks and discussions.
\end{ackn}

\bigskip

\bibliographystyle{amsplain}
\bibliography{HarmSoliton}

\bigskip
\bigskip

\end{document}